\documentclass[a4paper,12pt
]{article}
\usepackage{mathpazo}
\usepackage{amsmath}
\usepackage{amssymb}
\usepackage{latexsym}
\newtheorem{thm}{Theorem}[section]
\newtheorem{lem}[thm]{Lemma}
\newtheorem{prop}[thm]{Proposition}

\newtheorem{cor}[thm]{Corollary}
\newtheorem{rem}[thm]{Remark}
\newtheorem{df}[thm]{Definition}
\newtheorem{ex}[thm]{Example}

\newcommand{\R}{{\Bbb R}}
\newcommand{\N}{{\Bbb N}}
\newcommand{\C}{{\Bbb C}}
\newcommand{\Z}{{\Bbb Z}}

\newcommand{\F}{{\cal F}}

\newcommand{\La}{{\cal L}}

\newcommand{\MF}{{\cal S}}
\newcommand{\IE}{{\it i.e.}, }
\newcommand{\EG}{{\it e.g.}, }
\newcommand{\RESP}[1]{[{\it resp.} {#1}]}

\newcommand{\E}{E_\omega}
\newcommand{\nil}[1]{\mathit{Nil}^3({#1})}
\newcommand{\Or}{{\mathrm O}}
\newcommand{\cF}{\check{\bf F}}
\newcommand{\cC}{\check{\C}}
\newcommand{\olrho}{\overline{\rho}}
\newcommand{\ol}{\overline}
\newcommand{\WP}{P_{\bf w}^2}
\newcommand{\Lam}{F}
\newcommand{\AF}{{\mathbf{F}}}

\newlength{\itemwidth}
\begin{document}
\hspace*{\fill} (3rd revised version)
\begin{center}
{\large 
{\bf \LARGE 
Leafwise Symplectic Structures 
\vspace{10pt} \\
on Lawson's Foliation
} \vspace{20pt} \\
Yoshihiko MITSUMATSU\footnote{
Supported 
by Grant-in-Aid for 
Scientific Research (B) 22340015.
  
{\it Date}: October 31, 2013.   
 {\copyright}  2013 Y. Mitsumatsu.  

}
(Chuo University, Tokyo)
}\vspace{5pt}
\end{center}
\newlength{\abst}
\setlength{\abst}{120mm}

\begin{center}
\begin{minipage}{\abst}
\indent 
 The aim of this paper is to show that 
Lawson's foliation on the 5-sphere admits 
a smooth leafwise symplectic structure. 
The main part of the construction is 
to show that the Fermat type cubic surface 
admits an end-periodic symplectic structure. 
The results is paraphrased that the 5-sphere 
admits a regular Poisson structure 
of symplectic dimension 4. 
\end{minipage}
\end{center}

\setcounter{section}{-1}
\section{\large Introduction}
\indent 
In this article we show the following. 
\vspace{10pt}
\\
{\bf Theorem A} (Theorem \ref{MainTheorem}) \quad 
Lawson's foliation on the 
5-sphere $S^5$ admits 
a smooth leafwise symplectic structure. 

More generally, the Milnor fibration associated with 
one of the three classes of simple elliptic singularities 
$\tilde{E}_6, \tilde{E}_7$, and $\tilde{E}_8$ 
in complex three variables can be 
modified into a leafwise symplectic foliation on $S^5$ 
of codimension 1. 
Lawson's one is associated with $\tilde{E}_6$.  
\vspace{10pt}

We can paraphrase this result 
into the following. 
\vspace{10pt}
\\
{\bf Corollary B} \quad 
Associated with each of the three classes 
$\tilde{E}_6, \tilde{E}_7$, or $\tilde{E}_8$  
of simple elliptic singularities in three variables, 
the 5-sphere $S^5$ admits 
regular Poisson structure of symplectic 
dimension 4. 
\vspace{10pt}

This work is motivated and inspired 
by the works (\cite{SV}, [MV1, MV2]) 
by Alberto Verjovsky and others 
in which they are discussing the existence of 
leafwise symplectic and complex structures 
on Lawson's foliation and on its 
slight modifications. 
The author is extremely grateful to 
Verjovsky for drawing his attentions to such 
interesting problems.  

H. B. Lawson, {\sc Jr.}\ constructed 
a smooth foliation of codimension one 
on $S^5$ (\cite{L}), 
which we now call  {\it Lawson's foliation}. 
It was achieved by a beautiful combination 
of the complex and differential topologies and 
was a kind of breakthrough 
in an early stage of the history of foliations.  
The foliation is composed of two components.  
One is a 
tubular neighbourhood of a 
3-dimensional nil-manifold 
and the other one is, away from the boundary, 
foliated by Fermat-type cubic complex surfaces. 
As the common boundary leaf, 
here appears one of Kodaira-Thurston's 4-dimensional 
nil-manifolds.  
As each Fermat cubic leaf is spiraling 
to this boundary leaf, its end is diffeomorphic 
to a cyclic covering of 
Kodaira-Thurston's nil-manifold. 
(See Section \ref{Lawson} for the detail.)

In order to introduce a leafwise symplectic structure 
(for a precise definition, 
see Section \ref{Kodaira-Thurston}), 
we 
need 
to find a symplectic structure on 
the Fermat cubic surface 
which (asymptotically) coincides on the end 
with that of the cyclic covering 
of the Kodaira-Thurston nil-manifold. 
However, 
the natural symplectic structure 
on the Fermat cubic 
surface as 
an affine 
surface 
is quite different 
from the periodic ones on its end, because it is 
`conic' and expanding, not at all periodic.  
This is the crucial point in our problem.  
Once we find an end-periodic 
symplectic structure on the Fermat cubic surface 
(Section~\ref{end-periodic}), 
in fact it is almost enough to construct a smooth leafwise 
symplectic structure (Section~\ref{symplectic}), 
because 
it is easy to see 
that a 
simple foliation 
on the tube component admits 
a leafwise symplectic structure 
(Section~\ref{Kodaira-Thurston}). 


This paper is organized as follows. Lawson's foliation is 
reviewed in Section~\ref{Lawson}. 
In Section~\ref{Kodaira-Thurston}, 
certain symplectic structures on 
Kodaira-Thurston's nil-manifold and its covering 
are presented. 
According to the choice of the symplectic structures 
on Kodaira-Thurston's nil-manifold, 
a leafwise symplectic structures on the tube component 
is constructed.  
This enables us to state our real problem 
precisely in Section~\ref{symplectic}.  
Assuming the results on symplectic structures 
on the Fermat cubic surface in the later sections, 
a construction of a leafwise symplectic structure 
on Lawson's foliation is given in this section. 
Then, the natural 
symplectic structure on the Fermat cubic surface 
is analyzed in Section~\ref{FermatCubic}.  
Based on this analysis, in Section\ref{end-periodic}, 
which is the essential part of the present article,  
the existence of an end-periodic symplectic structure 
on the Fermat cubic surface is shown. 
In Section 6 we remark that our construction holds 
almost verbatimly in two other cases 
of the simple elliptic hypersurface singularities.  
In the final section, some related topics 
concerning the method in this paper are discussed. 

%
%
%
%
%
%
%
%
%
%
%
%
%
%
%
%
%
%
%
%

\section{\large Review of Lawson's Foliation}
\label{Lawson}
First we review the structure 
of Lawson's foliation $\La$. 
For those who are familiar with the materials 
it is enough to check 
our notations, 
which are in fact quite 
different from \cite{SV} and even from \cite{L}. 

Let us take a Fermat type homogeneous cubic polynomial 
$f(Z_0, Z_1, Z_2)=Z_0^3+Z_1^3+Z_2^3$ in three variables 
$Z_0, Z_1$, and $Z_2$. The complex surface ${\bf F}(w)=
\{(Z_0, Z_1, Z_2)\in\C^3 ; f(Z_0, Z_1, Z_2)=w\}$ for 
a complex value $w$ is non-singular if $w\ne0$ and 
${\bf F}_0$ has the unique singularity at the origin. 
The scalar multiplication 
$c\cdot (Z_0, Z_1, Z_2) = (cZ_0, cZ_1, cZ_2)$ 
by $c\in \C$ 
maps ${\bf F}_w$ to ${\bf F}_{c^3w}$. 
Hence ${\bf F}_0$ is preserved by such homotheties 
and for $w\ne0$ ${\bf F}_w$ is preserved iff $c^3=1$.

Now we put $ 
\tilde{F}_\theta=\bigcup_{\arg w = \theta}{\bf F}_w$ and 
$F_\theta=\tilde{F}_\theta \cap S^5$ where $S^5$ 
denotes 
the unit sphere in $\C^3$. 
Also, we put $N=\mathit{Nil}^3(-3) = {\bf F}_0 \cap S^5$. 
Let $p$ and $h$ denote the projection 
$p:\cC^3 \to S^5$ 
and the Hopf fibration $h : S^5 \to \C P^2$. 
Here $\cC^3$ denotes $\C^3 \setminus \{\Or\}$. 
Sometimes $h$ also denotes the composition 
$h\circ p : \cC^3 \to \C P^2$. 
Let also $H(t)$ $(t\in \R/2\pi\Z)$ denote 
the {\it Hopf flow} 
obtained by  scalar multiplication by 
$e^{it}$, whose orbits are the Hopf fibres.   
$E_\omega=\{[z_0:z_1:z_2];z_0+z_1+z_2=0\}=h({\bf F}_0)$ 
is an elliptic curve in $\C P^2$ with modulus 
$ 
\omega=\frac{-1 + \sqrt{-3}}{2}$.  
The Hopf fibration restricts to
$N \to E_\omega$, which is 
an $S^1$-bundle with $c_1=-3$. 
Also 
put 
$\MF_r=\cup_{\theta \in \R/2\pi\Z}{\bf F}_{re^{i\theta}}
= \{
Z=(Z_0, Z_1, Z_2)\in\C^3 ; \vert f(Z)\vert=r
\}
$ for $r>0$.  $\MF_r$ is also preserved by the Hopf flow.  

The following facts are also easy to see, 
while they are listed as Proposition 
for the sake of later use.  

\begin{prop} 
{\rm  1) 
$f|_{S^5}$ has no critical points around 
$N$. 
\vspace{2pt}
\\
\hspace*{8pt}
\begin{minipage}{\itemwidth}{
\hspace*{-15pt}
2) 
$\arg\circ f|_{S^5\setminus N}:S^5\setminus N \to S^1$ 
has no critical points away from $N$ and is called 
the {\it Milnor fibration}.  
Each fibre is by $F_\theta$ ($\theta\in\R/2\pi\Z$).  
\vspace{1pt}
}
\end{minipage}
\\
\hspace*{8pt}
\begin{minipage}{\itemwidth}{
\hspace*{-15pt}
3) The projection 
$p|_{{\bf F}_w} : {\bf F}_w \to F_\theta$ (${\theta =\arg w}$) 
is a diffeomorphism for $w\ne 0$.  
\vspace{3pt}
}
\end{minipage}
\\
\hspace*{8pt}
\begin{minipage}{\itemwidth}{
\hspace*{-15pt}
4) 
The Hopf fibration 
$h|_{F_\theta} : F_\theta \to {\C}P^2 \setminus \E$ 
restricted to $F_\theta$ 
is a three fold regular covering, and so is 
${\bf F}_w \to {\C}P^2 \setminus \E$ for $w\ne 0$. 
\vspace{1pt}
}
\end{minipage}
\\
\hspace*{8pt}
\begin{minipage}{\itemwidth}{
\hspace*{-15pt}
5) The normal bundle to $N\hookrightarrow S^5$ is 
trivialized by the value of $f$.  
\vspace{3pt}
}
\end{minipage}
\\
\hspace*{8pt}
\begin{minipage}{\itemwidth}{
\hspace*{-15pt}
6) The projection 
$p|_{\MF_r} : {\MF}_r \to S^5 \setminus K$ 
is a diffeomorphism for $r > 0$ 
and is equivariant with respect to the Hopf flow. 
$\MF_r$ fibers over the circle with fibres 
${\bf F}_{re^{i\theta}}$ ($\theta\in\R/2\pi\Z$)
and the fibration structure coincides with that of the 
Milnor fibration through $p|_{\MF_r}$.  
\vspace{3pt}
}
\end{minipage}
\\
\hspace*{8pt}
\begin{minipage}{\itemwidth}{
\hspace*{-15pt}
7) $H(2\pi/3)$ gives the natural monodromy 
of the Milnor fibration. 
}
\end{minipage}

}
\end{prop}
%
%
%
%
%
%
%


Let $W_r$ denote the tubular neighbourhood 
$W_r= \{Z=(Z_0, Z_1, Z_2)\in S^5 ; 
|f(Z)|\leq r\}$ 
of $N \subset S^5$ for $0<r$.  
Take $\varepsilon$ ($0<\varepsilon \ll 1$) so small 
that $f\vert_{W_\varepsilon}$ has no critical points. 
(Later we will take $\varepsilon$ again smaller 
for some reason.  )
Then let $U_\varepsilon=h(W_\varepsilon)\subset {\C}P^2$ 
be a tubular neighbourhood 
of $\E \subset {\C}P^2$.  
$W_\varepsilon$ is invariant under the Hopf flow. 
We choose further smaller constants 
$r_0$ and $r_*$ 
satisfying 
$0<r_0<r_*<\varepsilon$ and take 
$W_*=W_{r_*}$, $W=W_{r_0}$ 
and $U_*=U_{r_*}$, $U=U_{r_0}$. 
We decompose $S^5$ into ${W}$ and 
$C=S^5\setminus \mathrm{Int}\,W$, 
each of which are called 
the {\it tube} component and 
the {\it Fermat cubic} component. 
The statement 5) in the above proposition 
tells that $W_r$ is 
diffeomorphic to the product $N\times D^2_r$, 
while $\E\hookrightarrow {\C}P^2$ is twisted 
because $[\E]^2=9$. 
Here 
$D^2_r$ 
denotes the disk of radius $r$ in $\C$. 

The common boundary $\partial W= \partial C$ is 
diffeomorphic to $N\times S^1$, which is one of 
Kodaira-Thurston's 4-dimensional nil-manifolds 
and is well-known to be non-K\"ahler 
because $b_1=3$.  
It admits symplectic structures as well as 
complex structures 
but both are never compatible. 

As the two components $W$ and $C$ 
are fibering over the circle, 
the following lemma 
(a standard process of {\it turbulization}) 
is enough 
in order 
only to obtain a smooth foliation. 
However, to put leafwise symplectic structures, 
it is helpful to describe the foliation and 
the turbulization in more detail.   

\begin{lem}{\rm (\cite{L}, Lemma 1) \quad  
Let $M$ be a compact smooth manifold with boundary 
$\partial M$ 
and  
$\varphi: M \to S^1$ be 
a smooth submersion to the circle.  
Accordingly 
so is $\varphi\vert_{\partial M} : \partial M \to S^1$. 
Then, there exists a smooth foliation 
of codimension one 
for which the boundary $\partial M$ 
is the unique compact leaf, 
other leaves are diffeomorphic to 
the interior of the fibres, 
and the holonomy of the compact 
leaf is trivial as a $C^\infty$-jet.  
If we have two such submersions 
$\varphi_i: M_i \to S^1$ ($i=1,2$) 
with diffeomorphic boundaries 
$\partial M_1 \cong \partial M_2$, 
then on the closed manifold 
$M_1\cup_{\partial M_1 = \partial M_2} M_2$, 
by gluing them  
we obtain a smooth foliation of codimension one. 
}
\end{lem}

Let us formulate the turbulization process 
more explicitly. 
Firt take small positive constants 
$0<r_0<r_1<r_2<r_*$ and 
smooth functions $g(r)$ and $h(r)$ on $\R_+$ 
satisfying the following conditions.  
$$\begin{array}{rlrl}
g \equiv 0 
&\quad (r \leq r_0) , 
&
\qquad \qquad 
h \equiv 1
& 
\quad (r \leq r_1) ,  
\\
g = -\frac{\log \pi}{3}r 
&
\quad (r_1 \leq r <r_*) , 
&
\qquad \qquad 
h \equiv 0
& 
\quad (r_2 \leq r <r_*) ,
\\
g' < 0 
&
\quad (r_0 < r < r_*) , 
&
\qquad \qquad 
h' < 0
& 
\quad (r_1 < r < r_2) .  
\end{array}
$$
Then take a smooth non-singular vector field  
$X = g\frac{\partial}{\partial r} 
+ h\frac{\partial}{\partial \theta}$  
on the punctured plane 
$\R^2_+ = \R_+ \times \R/2\pi\Z$ where 
$(r, \theta)$ denote the polar coordinates.  
The integral curves of $X$ define 
a smooth foliation $\F_T$ 
on 
$\C\setminus \{\Or\}\cong\R_+ \times \R/{2\pi\Z}$. 
The constant $\frac{\log \pi}{3}$ has no significance 
at this stage.

Now the turbulization on the side 
of the Fermat cubic component 
is described as follows.  
The foliation 
$\tilde\La=\{\Lam_\theta\}$ on $S^5\setminus N$ 
by the Milnor fibres and the pull-back foliation 
$f^{-1}\F_T$ coincide with each other 
on $W_*\setminus W_{r_2}$. 
On the Fermat cubic component 
$C= S^5\setminus \mathop{\mathrm{Int}}W$, 
Lawson's foliation $\La\vert_C$ is obtained 
as 
$\La\vert_{S^5\setminus W_{r_2}}
=\tilde{\La}\vert_{S^5\setminus W_{r_2}}$ 
and  
$\La\vert_{W_*\setminus W}
=f^{-1}\F_T\vert_{W_*\setminus W}$.  
Let $L_\theta$ denote one of the resulting 
leaves which contains 
$\Lam_\theta\setminus W_{r_2}$.  
$L_\theta$ is diffeomorphic to $\Lam_\theta$ and 
only the embedding of the product end 
$N\times \{r\cdot e^{i\theta}\}$ is modified 
by the turbulization procedure. We will 
fix an identification 
of $L_\theta$ with ${\bf F}_{e^{i\theta}}$ 
in Section~\ref{symplectic}.

On the tube component, we can also formulate 
the turbulization using $\F_T$ in a similar way 
to the above, 
however,  
since $\La\vert_W$ 
can also described using a well-known foliation, 
the ``Reeb component'',  
on $S^1\times D^2$, 
here we adopt such a description. 
The tube component $W$ is diffeomorphic to 
$N\times D^2$ and $N=\nil{-3}$ is an $S^1$-bundle 
over the elliptic curve $\E$. 
We take a smooth coordinate $(x, y)$ for $\E$ 
where $x, y \in S^1=\R/{2\pi\Z}$. 
Then the projection from $\E$ to $S^1\ni x$ 
gives rise to a fibration of $N$ over $S^1$ with 
fibre $T^2$ and the monodromy 
${1\, -3}\choose{0\,\,\,\,\, 1}$. 
Therefore the tube component $W$ fibers 
over the solid torus $S^1\times D^2$ 
with the fibre $T^2$ 
and the monodoromy 
${1\, -3}\choose{0\,\,\,\,\, 1}$, 
namely, $W=\R\times D^2 \times T^2/\sim$ 
where 
$(x+2\pi, P, {y\choose z})
\sim
(x, P, {{1\, -3}\choose{0\,\,\,\,\, 1}} 
{y\choose z})$. 

As the tube component part of Lawson's foliation $\La$, 
we can take the pull-back of the standard Reeb component 
$\F_R$ 
on $S^1 \times D^2$ to $W$. 
Thus we obtain Lawson's foliation $\La$ on $S^5$, 
whose boundary is a unique compact leaf and 
is diffeomorphic 
to Kodaira-Thurston's 
nil-manifold .  
\begin{rem}\label{nil}
{\rm \quad The 3-dimensional nil-manifold 
$\nil{c_1}$ is often presented as the quotient 
$\nil{c_1}=\Gamma(c_1)\setminus H$ 
of the 3-dimensional Heisenberg group $H$ 
by its lattice $\Gamma(c_1)$, which are defined as 
\vspace{-5pt}
$$H=\{\left(\begin{array}{ccc}
1&\ol{x}&\ol{z}
\\
0&1&\ol{y}
\\
0&0&1
\end{array}
\right)
\,;\,\ol{x},\ol{y}, \ol{z} \in \R\} \supset  
\Gamma(c_1)=\{\left(\begin{array}{ccc}
1&\ol{x}&\ol{z}
\\
0&1&\ol{y}
\\
0&0&1
\end{array}
\right)
\,;\,\ol{x},\ol{y}, c_1\ol{z} \in \Z\}.
$$ 
In the case $c_1<0$, $\ol{z}$ must be understood 
to have opposite sign. 
In this coordinate on $H$ 
take $\frac{\partial}{\partial \ol{x}}$, 
$\frac{\partial}{\partial \ol{y}}$, 
and 
$\frac{\partial}{\partial \ol{z}}$ 
at the unit element, and then extend them 
to be $X$, $Y$, and $Z$ as left 
invariant vector fields.  
Let $d\ol{x}$, $d\ol{y}$, and $\ol{\zeta}$ be 
the dual basis for the invariant 1-forms, which 
satisfies  $d\ol{\zeta}=d\ol{x}\wedge d\ol{y}$.  
On our 
$N=\nil{-3}$ we have $x=2\pi\ol{x}$, $y=2\pi\ol{y}$, 
$z= 2\pi c_1 \ol{z}$, 
and $\zeta= 2\pi c_1 \ol{\zeta}$. 
}
\end{rem}
\begin{rem}\label{N'}
{\rm \quad 
As the normal bundle to $\E \hookrightarrow {\C}P^2$ 
has $c_1=9$, the boundary $\partial U$ is isomorphic to 
$\nil{9}$.  Here we are looking at $\partial U$ 
from the interior of $U$.  
However, for the later purpose, it is more convenient 
to give the opposite orientation, because 
to $N$ we gave the orientation as the boundary of 
${\bf F}_0\cap\{\rho \leq R\}$ as well as the end of ${\bf F}_1$.  
Therefore let $N'$ denote 
$\partial ({\C}P^2 \setminus U)$, 
which is isomorphic to $\nil{-9}$. 
}
\end{rem}

\section{\large 
Symplectic Forms on the Kodaira-Thurston 
Nil-Manifold and ${\bf F}_0$} 
\label{Kodaira-Thurston}
In this section, 
we describe natural symplectic forms 
on Kodaira-Thurston's 4-dimensional nil-manifold 
and show that the tube component admits 
a smooth leafwise symplectic structure 
which is {\it tame}  
around the boundary. 

\begin{df}{\rm \quad 
A smooth 
{\it leafwise symplectic structure} (or {\it form}) on 
a smooth foliated manifold $(M, \F)$ 
is a smooth leafwise closed 2-form $\beta$ 
which is non-degenerate on 
each leaves. 

More precisely, first, 
$\beta$ is a smooth section to the smooth vector 
bundle $\bigwedge^2T^*\F$. 
For smooth sections to  $\bigwedge^*T^*\F$ naturally 
the exterior differential in each leaves is defined. 
This exterior differential 
is often  denoted by $d_\F$. 
$\beta$ is closed in this sense 
and 
is non-degenerate in each leaves, 
namely, $d_\F\beta=0$ holds and 
 $\beta^{\dim\F/2}$ defines 
a volume form on each leaves.  

The existence of such $\beta$ is equivalent to 
that of a smooth 2-form $\tilde\beta$ on $M$ 
whose restriction to each leaf 
is a symplectic form of the leaf. 
It should be remarked that 
$\tilde\beta$ may not be closed as a 2-form on $M$.  
Furthermore, quite often we mix the two formulations 
and do not make a clear distinction. }
\end{df}
\begin{df}\label{tame-def}{\rm \quad 
Let $(M, \F)$ be a smooth foliated manifold 
with a boundary compact leaf $\partial M$ and 
a leafwise symplectic form $\beta$. 
$(M, \F, \beta)$ is {\it tame} around the boundary 
if the triple satisfies the following condition. 
We also simply say that $\beta$ is tame. 
\begin{enumerate}
\item[(1)] The (one-sided) holonomy 
of the boundary leaf 
is trivial as $C^\infty$-jet. 
\item[(2)] There exists a collar neighbourhood 
$V\cong [0,\epsilon)\times \partial M$ 
of the boundary $\partial M$ with the projection 
$\mathit{Pr}: [0,\epsilon)\times \partial M \to \partial M$ 
for which 
$\beta \vert_V$ coincides 
with the restriction to the leaves 
of the pull-back 
$\mathit{Pr}^*(\beta\vert_{\partial M})$. 
\end{enumerate}
}
\end{df}

\begin{cor}\label{tame-cor}
{\rm \quad 
Let $(M_i, \F_i, \beta_i)$ ($i=0,1$) be  two 
foliated manifolds with leafwise symplectic structures.   
Assume that both are tame around their boundaries 
and there exists a symplectomorphism 
$\varphi: (\partial M_1, \beta_1\vert_{\partial M_1}) 
\to 
(\partial M_2, \beta_2\vert_{\partial M_2})$ 
between their 
boundaries. 
Then gluing by $\varphi$ yields 
a smooth foliated manifold 
$(M=M_1\cup_\varphi M_2, \F, \beta)$ 
with a smooth leafwise symplectic structure. }
\end{cor}

We adopt this corollary for our construction.  
The existence of tame leafwise symplectic structures 
on the Fermat cubic component is discussed 
in Section 3, 4, and 5. 
In this section we show the existence of tame ones on the 
tube component. 

For the 3-dimensional Reeb component $(S^1\times D^2, \F_R)$,  
as the leaves are 2-dimensional,   
it is easy to show that 
there exists a tame leafwise symplectic structure. 
For any area form of the boundary, 
extend it to a collar 
neighbourhood so as to satisfy the tameness condition, 
and then further extend it to a leafwise 2-form 
on the whole component 
so that on each leaf it gives an area form. 
In this construction we can start 
with the standard area form 
$dx \wedge d\theta$ of the boundary 
under the coordinates 
$(x, r, \theta)$ for 
$S^1\times D^2=\{(x, r, \theta)$ $;$ 
$x\in S^1=\R/{2\pi\Z},$ 
$0\leq r \leq r_0,$ 
$\theta\in S^1=\R/{2\pi\Z}\}$.  
The resulting tame leafwise symplectic form on 
 $(S^1\times D^2, \F_R)$ is denoted by $\beta_R$. 

Now let us construct tame leafwise symplectic forms 
on the tube component, 
using the description of the tube component 
in the end of the previsous section. 
%
Let $\zeta$ denote 
the standard connection 1-form for the Hopf fibration 
$h : S^5 \to {\C}P^2$.  $\zeta$ coincides 
with the standard contact form 
$\zeta = \sum_{j=1}^{3}(x_jdy_j-y_jdx_j)$ 
on $S^5$.  
On each fibre (with an arbitrary reference point) 
$\zeta$ defines an identification with 
$S^1=\R/{2\pi\Z}$ and 
the resulting coordinate 
is denoted by $z$ in the previous section. 
Once $\zeta$ is restricted to $N=F \cap S^5$ 
it is denoted by $\zeta_N$.

The tube component $W$ admits 
a flat bundle structure 
$T^2\hookrightarrow W=N\times{D^2} 
\stackrel{\pi_{\!R}}{\to} S^1\times{D^2}$ 
with the monodromy 
${{1\, -3}\choose{0\,\,\,\,\, 1}}$.  
On $N$ we have 
$d\zeta_N=-(-3)2\pi\frac{dx}{2\pi}
\wedge\frac{dy}{2\pi}
=\frac{3}{2\pi}dx\wedge{dy}
$ 
and hence 
$dy\wedge\zeta_N$  
is a closed 2-form which restricts 
to a holonomy invariant area form 
on each fibre $\cong T^2$.

As Lawson's foliation on the tube component $W$ is 
given as ${\pi_{\!R}}^{-1}\F_R$, 
using the pull-back ${\pi_{\!R}}^*\beta_R$ 
we obtain a tame leafwise symplectic form 
$\beta_{W, \lambda,\mu}= 
\lambda\,{{\pi_{\!R}}^*\beta_R} 
+ \mu\,{dy\wedge\zeta_N}$ 
for non-zero constants $\lambda$ and $\mu$.  
The restriction of $\beta_{W, \lambda,\mu}$ 
to the boundary $\partial W$ 
is presented  as 
$\lambda\,{d\theta\wedge{dx}} 
+ \mu\,{dy\wedge\zeta_N}$.  
Also it is easy to see that the foliation $\La\vert_W$ 
(in fact $\La$ itself)
and the leafwise symplectic form  
$\beta_{W, \lambda,\mu}$ 
is invariant under the Hopf flow $H(t)$.   
We have established the following. 

\begin{prop}\label{TubeComp}{\rm \quad 
On the tube component Lawson's foliation 
$\La\vert_W$ admits a tame leafwise symplectic form 
$\beta_{W, \lambda,\mu}= 
\lambda\,{{{\pi_{\!R}}^{}}^*\beta_R} 
+ \mu\,{dy\wedge\zeta_N}$ 
for constants $\lambda\ne0$ and $\mu\ne0$,  
which is invariant under the Hopf flow. 
It restricts to 
$\beta_{K, \lambda,\mu}= 
\lambda\,{d\theta\wedge{dx}} 
+ \mu\,{dy\wedge\zeta_N}$
on the boundary leaf.  
}
\end{prop}

In order to make a better correspondence with 
the Fermat cubic component, 
we introduce a new coordinate variable $\tau$ 
which replaces $\theta$ 
in the presentation of these symplectic forms 
on Kodaira-Thurston's nil-manifold.  
$\theta$ is reserved for the argument 
of the value of $f$.  
Then consider a new relation 
$\tau=2\log \rho$ 
between the radius 
$\rho=\sqrt{
\vert Z_0\vert^2+\vert Z_1\vert^2
+\vert Z_2\vert^2}$ on $\C^3$,   
which fits into the following picture.  
Since 
our Kodaira-Thurston nil-manifold $K$ 
is presented as 
$K=\partial W = N \times S^1$ 
and also 
the polar coordinate $\R_+ \times S^5$ 
for $\cC^3$ restricts 
to 
$\cF_0={\bf F}_0\setminus\{\Or\}
= \R_+ \times N$, 
we identify $K$ with $\cF_0/\!_\sim$ 
where 
$P \sim Q$ 
for $P$ and $Q \in \cF_0$ 
iff 
$Q=e^{n\pi}P$ for some $n\in\Z$. 
\begin{rem}
{\rm \quad Instead of $e^{\pi}$, 
we can take any complex number 
$c\in\C$ with $\vert c\vert>1$. Then according to $c$ 
the complex structure 
induced on $K$ from ${\bf F}_0$ varies. However, they are all 
diffeomorphic to each other and in this article 
we are indifferent to the complex structures and 
rather interested 
in the symplectic structures of $K$.  
}
\end{rem}

On the Fermat cubic component $C$, 
any of the interior leaf 
is diffeomorphic to the Fermat cubic surface ${\bf F}_1$, 
and is spiraling to the boundary leaf $K$. 
The way in which an interior leaf approaches to $K$ 
is topologically the same as ${\bf F}_1$ approaches 
to ${\bf F}_0$ after divided by the scalar multiplication 
by $e^{n\pi}$ for $\forall n \in \Z$. 
Looking at $\cF_0$ and ${\bf F}_{e^{i\theta}}$'s 
in this way is used in later sections.   
It is also helpful in introducing 
leafwise complex structures 
on the Fermat cubic component.  
In the rest of this section, 
we take a closer look at   
$K=\cF_0/\!_\sim$ and 
$\cF_0 = \R_+ \times N$.   

From $K$ the infinite cyclic covering pulls 
$\beta_{K, \lambda,\mu}$ back 
to a periodic symplectic form   
$\beta_{0, \lambda,\mu}= 
\lambda\,{d\tau\wedge{dx}} 
+ \mu\,{dy\wedge\zeta_N}
= \lambda\,{\frac{2}{\rho}d\rho\wedge dx} 
+ \mu\,{dy\wedge\zeta_N}
$
on $\cF_0$.  
On the other hand, $\cF_0$ 
inherits a natural symplectic structure 
from 
$(\C^3, \beta^*=2\sum_{j=0}^2dx_j\wedge dy_j)$.  
One of the standard Liouville forms 
(the primitive of symplectic form) 
$\lambda^*=\sum_{j=0}^2(x_jdy_j - y_j dx_j)$ 
is presented as 
$\rho^2\zeta_N$ in the polar coordinate. 
Replacing $\rho^2$ with $\olrho$, 
we see that $(\cF_0, \beta_0=d(\olrho \zeta_N))$ is 
the symplectization of the contact manifold 
$(N, \zeta_N)$.  
Replacing $\rho$ with 
$\tau= \log\olrho=2\log\rho \in \R$, 
we have an identification of 
$\cF_0$ with $\R \times N \ni (\tau, P)$.  
We call this the {\it product}   
coordinate. 
%
%
%
%
%
%
%
%
%
%
%
%
%
%
%
%
%
%
%
%
%
%
%
%
%
%
%
%
%
%
%
%
%
\section{\large 
Leafwise Symplectic Structure 
on 
the Fermat Cubic Component}
\label{symplectic}
\subsection{Main Theorem}
\begin{thm}\label{MainTheorem}{\rm\quad
Lawson's foliation $\La$ on $S^5$ admits 
a smooth leafwise 
symplectic form. 
}
\end{thm}
This is the main theorem of the present article.  
Proposition~\ref{TubeComp} and Corollary 2.3 
imply that this is a direct consequence 
of the following proposition, which we prove 
in this section 
assuming Corollary~\ref{end-periodicCor}.  
\begin{prop}\label{MainProposition}
{\rm\quad
For a sufficiently small constant 
$0<\mu\ll1$ 
and a sufficiently large constant 
$\lambda\gg1$ 
Lawson's foliation restricted 
to the Fermat cubic component 
admits a tame symplectic form which 
restricts to 
$\beta_{K, \lambda,\mu}= 
\lambda\,{d\tau\wedge{dx}} + \mu\,{dy\wedge\zeta_N}$
on the boundary leaf $K$.  
}
\end{prop}
\subsection{Coordinates}

As a preparation, 
we start with a more precise description 
of tubular neighbourhoods of $N$ in $S^5$ 
and those of $\cF_0$ in $\cC^3$.  
For the tubular neighbourhood $W_\varepsilon$ 
of $N$ in $S^5$, 
we describe an identification 
of $W_\varepsilon$ with $D^2_\varepsilon\times N$. 
As $f:W_\varepsilon\to D^2_\varepsilon\subset \C$ 
defines the first projection, 
it is enought to define 
the projection to $N$. 

In Section 1 we chose $\varepsilon$ so small that 
$W_\varepsilon$ avoids the critical points 
of $f\vert_{W_\varepsilon}$.   
Here we take $\varepsilon$ again smaller in the following sense. 
Let $\nu N$ denote the normal bundle to 
to 
$N\subset S^5$.  
With respect to the standard metric induced from $\C^3$, 
there exists a radius $\delta>0$ for $\nu N$ such that 
the exponential map $\nu N \to S^5$ is diffeomorphic 
up to the radius $\delta$. 
We choose $\varepsilon > 0$ so that 
$W_\varepsilon$ is included in 
this diffeomorphic image.  
Thus the inverse of the exponential map 
defines the projection from $W_\varepsilon$ to $N$. 
Simply, this projection asigns to any point in $W_\varepsilon$ 
its nearest point in $N$. 
As the Hopf action is isometric and leaves $N$ invariant, 
this projection is equivariant with respect to the Hopf action.



We do the same on each 
sphere $S^5(\rho)$ of the radius $\rho \geq 1$.    
Remark that once $\varepsilon$ has been chosen 
small enough for $\rho=1$, so is it 
for any $\rho \geq 1$. 

Collecting these identifications with respect to 
each $\rho\geq 1$, 
it also defines an identification of the 
tubular neighbourhood 
$\cC^3
\cap\{\rho \geq 1\}
\cap\{\vert f\vert <\varepsilon\}$ 
of $\cF_0\cap\{\rho \geq 1\}$ 
with  
$(\cF_0\cap\{\rho \geq 1\})
\times D^2_\varepsilon$ 
and with $[1,\infty) \times N \times D^2_\varepsilon$. 
The end 
${\bf F}_w\cap\{ \rho \geq 1\}$  
of the Fermat cubic surface 
${\bf F}_w$ for $\vert w \vert \leq \varepsilon$ is 
exactly the graph 
of the constant function {$\equiv w$}.  
From this we see 
that these ends are diffeomorphic to 
$ [1, \infty)  \times N$, 
which is called the {\it product coordinate} of the end.   
Later, $\tau=2\log \rho$ replaces $\rho$, 
where the identification is 
denoted by 
$[0, \infty)\times N
\ni (\tau, P)$ 
and is also called the product coordinate.

\begin{rem}
{\rm \quad
We have remarked that the identification 
of the end of ${\mathbf F}_w$ with that of $\cF_0$ is 
equivariant with respect to the Hopf flow.  
Therefore the product coordinate is also 
equivariant under the Hopf flow $H(t)$ 
if we interpret the action of $H(t)$ 
on 
$
\R_+\times N \times D^2_\varepsilon$ 
is 
$(\rho, P, z) \mapsto 
(\rho, H(t)\cdot P, e^{i3t}z])$  
where $H(t)\cdot P$ is the restriction of the Hopf 
flow on $N\subset S^5$ and is nothing but 
the fibrewise multiplication by $e^{it}$ 
on each fibre of $S^1\hookrightarrow N \to \E$. 
%
Similarly, on $W_\varepsilon\setminus N \,
(\cong 
N \times D^2_\varepsilon)
\cong
N \times (0, \varepsilon) \times S^1
$, 
the Hopf action is indicated as 
$
H(t)(P, r, \theta)=(H(t)\cdot P, r, \theta + 3t)
$.  
}
\end{rem}
\begin{rem}\quad{\rm
For any $0<\epsilon\leq \varepsilon$,  
as $\mathbf{F}_{\epsilon}$ triply covers 
${\C}P^2\setminus E_\omega$ and 
the product coordinate is equivariant with respect to 
the monodromy $H(\frac{2\pi}{3}\Z)$, the product coordinate 
$\mathbf{F}_{\epsilon} \cong [0,\infty)\times N$ 
induces a product coordinate of the end of   
${\C}P^2\setminus E_\omega \cong [0,\infty)\times N'$, 
where 
$N'=-\partial U
\cong N/_{\Z/_3}\cong\mathit{Nil}^3(-9)$.  
}
\end{rem}
\subsection{Proof of main proposition}

First let us introduce a leafwise symplectic sturcture 
on the Milnor fibration 
$(S^5 \setminus N,\, \{\Lam_\theta\})$. 
Following 
Corollary \ref{end-periodicCor}, 
take and fix a symplectic form $\beta_{\lambda, \mu}$ 
on ${\bf F}_{\varepsilon'}$.  
(The constant ${\varepsilon'}$ is given in Theorem 
\ref{naturalThm} and Proposition \ref{symplemb}.)
Each leaf $\Lam_\theta$ is identified with ${\bf F}_{\varepsilon'}$ 
\vspace{-4pt}
through the projection and 
the Hopf actions; 
${\bf F}_{\varepsilon'}\,{\stackrel{p}{\longrightarrow}\,
}\Lam_0
\,
{\stackrel{H(t)}{\longrightarrow}}\,
\Lam_\theta$ for $t=\frac{\theta + 2k\pi}{3}$ 
(mod $2\pi\Z$) ($k=0,\, 1,\, 2$).  
Let $p_{\theta, k}:{\bf F}_{\varepsilon'} 
\to \Lam_\theta$ denote this 
identification.  
For the sake of continuity 
we can not decide which 
one to choose among three values of $t$. 
However as $\beta_{\lambda, \mu}$ is invariant 
under the action 
of $H(\frac{2\pi}{3}\Z)$ on ${\bf F}_{\varepsilon'}$, 
these indentifications induce 
a well-defined symplectic form on each $\Lam_\theta$ 
from  (${\bf F}_{\varepsilon'}, \, \beta_{\lambda, \mu}$), 
which gives rise to a smooth leafwise symplectic structure 
$\beta_{\tilde{\La}}$ on the Milnor fibration.    
\vspace{5pt}

Next we go back to the turbulizing process of obtaining 
the Fermat cubic component. 
Here we need a pointwise identification of each interior 
leaf $L_\theta$ of the Fermat cubic component with 
$\Lam_\theta$ and with $\Lam_0$. 
On $W_{r^*}\setminus N 
\cong N \times (D^2_{r^*}\setminus \{0\})$, 
take a vector field $\tilde X$ 
which is defined as $(0, X)$, where 
$X$ is the vector field defined on $D^2_{r^*}$ 
in Section \ref{Lawson} for the turbulization.  
Also take a vector field 
$\tilde R=(0, -\frac{\log \pi}{3}r\frac{\partial}{\partial r})$.

We identify  $L_0$ with $\Lam_0$ and 
thus with ${\bf F}_{\varepsilon'}$ as follows. 
The core part 
$L_0\setminus W_{r_2}$ 
is exactly identical with 
$\Lam_0\setminus W_{r_2}
$.   
For 
%
$t >-\frac{3}{\log \pi}(\log r_*-\log r_2)$,  
the point 
$\exp(t\tilde X)(P, r_2, 0)$ of $L_0$ 
is identified with the point 
$\exp(t\tilde R)(P, r_2, 0)$ of $\Lam_0$.  
Accordingly, these points are identified with 
the point $(\tau_2 + t, P)$ of ${\bf F}_1$ in the 
product coordinate, 
where 
$\tau_2=-\frac{2}{3}\log r_2$. 
Similarly 
$\exp(t\tilde X)(P, r_2, \theta)$ 
$\in$  
$L_\theta$ and 
$\exp(t\tilde R)(P, r_2, \theta) 
\in
\Lam_\theta$ are 
identified. 
Let 
$\mathit{lm}_\theta:\Lam_\theta \to L_\theta$ 
denote this identification. 
Remark that through these identifications 
a point $(\tau , P)$ in the end of ${\bf F}_{\varepsilon'}$ 
for large enough $\tau \gg 0$ 
is sent to a point 
$\mathit{lm}_0\circ p(\tau , P)
=(P, r(\tau), \tau + c_0)$ 
in $L_0 \cap W_{r_2}$ for some function $r(\tau)$ 
and some constant $c_0$.  
Also we have 
$\mathit{lm}_\theta\circ p_{\theta, k}(\tau , P)
=(H(\frac{\theta + 2k \pi}{3})P, 
r(\tau), \tau + c_0 + \theta)$.  
\vspace{5pt}

The leafwise symplectic form $\beta_{\tilde{\La}}$ 
on the Milnor fibration 
is thus transplanted on 
the interior of the Fermat cubic component $C$ 
of Lawson's foliation $\La$ 
to be a leafwise symplectic form 
$\beta_C$.  
What remains to prove is that 
$\beta_C\vert_{W_{r_1} \setminus \ol{W}}$  
coincides with 
$\mathit{Pr}^*\beta_{K, \lambda, \mu}$ where 
$\mathit{Pr}$ denotes the projection 
of the end 
$W_{r_1} \setminus \ol W \cong N\times (r_0, r_1) \times S^1$ 
of $C$ to the boundary 
$\partial C \cong N\times \{r_0\} \times S^1 
\cong N\times S^1$(${}=K$).  
Then we obtain a tame symplectic form on $C$ with 
the restriction $\beta_{K, \lambda, \mu}$ to the boundary.  
By Corollary~\ref{tame-cor} the proof will be completed.  
\vspace{5pt}

From the above preparations, 
we see  that the composition of the maps 
$\mathit{Pr}\circ \mathit{lm}_\theta\circ p_{\theta, k}
:$ [the end of ${\bf F}_{\varepsilon'}$] $= (T, \infty)\times N 
\to 
\partial C =N\times S^1$ 
sends the points as $(\tau , P)\mapsto 
(H(\frac{\theta+2k\pi}{3})P, \tau + c_0 + \theta))$ 
for some $T \gg 0$.  
As is mentioned in Corollary~\ref{end-periodicCor},  
$\beta_{\lambda, \mu}\vert_{(T, \infty)\times N}$ is 
invariant under the Hopf flow and the $\tau$-translation, 
for any $\theta$ and $k\in \Z$,   
$(\mathit{Pr}\circ \mathit{lm}_\theta\circ p_{\theta, k})_*
\beta_{\lambda, \mu}\vert_{(T, \infty)\times N}$ coincide 
with each other and in fact with  
$\beta_{K, \lambda, \mu}$.   
This is what we had to prove.  
\hspace*{\fill} $\square$\ref{MainProposition}.

\section{\large Natural Symplectic Structure 
on Fermat Cubic Surface}
\label{FermatCubic}%

\subsection{Statements and notations}
In order to construct an end-periodic symplectic structure 
on $\AF_{\varepsilon'}$, we need to know fairy exactly a 
symplectic stucture that  $\AF_{\varepsilon'}$ a priori 
admits. It is not very difficult to show that 
the natural symplectic structure which 
$\AF_{\varepsilon'}$ inherits from $\C^3$ 
coincides with that of  $\cF_0$ on the product ends 
up to an isotopy inside  $\AF_{\varepsilon'}$. 
However, for our purpose the following result is enough, 
which is much easier to prove.

\begin{thm}\label{naturalThm}
{\rm\quad
For a sufficiently small $\varepsilon'>0$, 
there exists a symplectic form $\beta_1$ 
on the Fermat cubic surface $\AF_{\varepsilon'}$ 
which satisfies the following 
properties.    \vspace{-5pt}
\begin{enumerate}
\item[(1)] On the end 
$\AF_{\varepsilon'} \cap \{\tau\geq 2\pi  \}$ 
in the product coordinate $[0, \infty)\times N$,  
$\beta_1\vert_{[2\pi, \infty)\times N} 
= d(e^\tau\zeta_N)$.  
 \vspace{-5pt}
\item[(2)]  $\beta_1$ is invariant under 
the Hopf action 
$H(t)$ for $t\in \frac{2\pi}{3}\Z$.  
 \vspace{-5pt}
\end{enumerate}
}
\end{thm}
\begin{cor}\label{naturalCor}
{\rm\quad
For the same $\varepsilon'>0$ as above, 
there exists a symplectic form $\beta'$ 
on ${\C}P^2\setminus\E$ 
whose restriction to the the product end 
satisfies 
$\beta'\vert_{[2\pi, \infty)\times N'} 
= d(e^\tau\zeta_{N'})$ 
with respect to the product coordinate 
$[0, \infty)\times N'$. 
}
\end{cor}
See Remark 3.4 for the product coordinate for the end of 
${\C}P^2\setminus\E$ 
and Remark~\ref{N'} for $N'$.  
We explain some more notations on ${\C}P^2\setminus \E$. 
The natural contact 1-form $\zeta_{N'}$ is obtained 
as $\ol\zeta$ in Remark~\ref{nil} and 
is also obtained as the quotient 
$(N', \zeta_{N'})=(N, \zeta_N)/_{\Z/_3}$. 
As  ${\C}P^2\setminus \E$  is regarded as 
the quotient of 
${\bf F}_1$ by the restricted Hopf action by $\Z/_3$, 
on its end the product coordinate 
$\cong(T, \infty)\times N'$ 
is also naturally 
induced from the product coordinate 
$(T, \infty)\times N$ for the end of ${\bf F}_1$ 
by simply regarding $N'=N/_{\Z/_3}$.

%
%
\subsection
{Re-embedding of Fermat cubic surfaces}
\label{Re-embedding}
We prove Theorem \ref{naturalThm}. 
On the tubular neighbourhood 
$\{\vert f\vert \leq \varepsilon\}\cap \{\tau \geq 0 \}$ 
of $\cF_0 \cap  \{\tau \geq 0 \}$  
in $\C^3 \cap \{\tau \geq 0 \}$, take the product coordinate 
$[0, \infty)\times N \times D^2_\varepsilon$. 
The end of an affine surface $\{\tau \geq 1\} \cap \AF_w$ 
is the graph of the constant function 
$c_w\equiv w:[0, \infty)\times N  \to D^2_\varepsilon$ in the 
product coordinate. 

Take a smooth function $\phi:[0,\infty)\to [0,1]$ 
satisfying the conditions 
$$\phi(\tau) \equiv 1 \quad \mathrm{for} 
\quad \tau \in [0,1]\qquad \mathrm{and} 
\qquad \phi(\tau) \equiv 0 \quad 
\mathrm{for}
\quad \tau \in [\pi,\infty)$$ 
and consider the graphs of the functions 
$$ \sigma_w : [0, \infty)\times N \to D^2_\varepsilon,
\quad \sigma_w(\tau, P)= \phi(\tau)w$$ 
for $w\in D^2_\varepsilon$. 
The part  $\{\tau\geq \pi\}$ of the graph 
coincides with 
$\cF_0\cap\{\tau\geq \pi\}$ and thus symplectic as submanifold 
of $(\C^3, \beta^*)$.  
(See the last paragraph of Section 2 for $\beta^*$. ) 
The family $\{\sigma_w\, ; \, w \in D^2_\varepsilon \}$ 
apparently depends smoothly on $w$ and the graphs converge to 
$\AF_0\cap\{\tau\geq 0\}$ when $w\to 0$.  
Therefore the following proposition follows easily 
from the compactness of $[0,2\pi]\times N$. 
\begin{prop}\label{symplemb}
{\rm\quad
There exists $0<\varepsilon'\leq \varepsilon$ such that 
the graph of $\sigma_w$ for $\vert w \vert \leq \varepsilon'$ 
is a symplectic submanifold of $(\C^3, \beta^*)$. 
}
\end{prop}

For example, if we take $w=\varepsilon'$, the affine surface 
$\AF_{\varepsilon'}$ has another smooth embedding into $\C^3$ 
which is modified from the original one only on the end 
$\{\tau\geq 0\}$ by the graph of $\sigma_{\varepsilon'}$ 
which is again, a symolectic submanifold and coincides with 
$\cF_0$ on the end. Because the new embedding is respecting 
the product coordinate, 
it is equivariant under the Hopf action $H(t)$ for 
$t\in \frac{2\pi}{3}\Z$. 
The core part $\AF_{\varepsilon'}\cap \{\tau\leq 1\}$ 
is unchanged, it is also invariant under 
$H(\frac{2\pi}{3}\Z)$. 
It is clear that $\beta^*$ is invariant under the Hopf flow 
and $\beta^*\vert_{\cF_0}=d(e^\tau\zeta_N)$. 
Therefore the induced symplectic form 
$\beta_1$ from  $(\C^3, \beta^*)$ by the new embedding 
is the desired one in Theorem \ref{naturalThm}. 
\hspace{\fill}$\square \ref{naturalThm}.$

\section{\large End-Periodic Symplectic Form 
on Fermat Cubic Surface}
\label{end-periodic}
Based on the preparations that we have done 
in the preceding sections,   
we prove the following results, 
which are the core part of this article,   
while the arguments are surprisingly simple and easy. 
We use the product coordinate 
$(\tau, P')\in 
(-\frac{2}{3}\log \varepsilon,  \infty)\times N'
\cong 
U \setminus \E 
\subset {\C}P^2 \setminus \E$.  
\begin{thm}\label{end-periodicThm}
{\rm\quad 
For a sufficiently small constant $0<\mu\ll 1$ and 
sufficiently large constants $T\gg2\pi$ and 
$\lambda\gg1$, 
there exists a symplectic form $\beta'_{\lambda, \mu}$ 
on ${\C}P^2 \setminus \E$  which restricts to 
$\lambda\,d\tau\wedge dx + \mu\,dy\wedge \zeta_{N'}$ 
on its end $(T, \infty) \times N'$. 
}
\end{thm}
\begin{cor}\label{end-periodicCor}
{\rm \quad For the same constants as above, 
there exists a symplectic form $\beta_{\lambda, \mu}$ 
on the Fermat cubic surface ${\bf F}_{\varepsilon'}$  
which restricts to 
$\lambda\,d\tau\wedge dx + \mu\,dy\wedge \zeta_{N'}$ 
on its end $(T, \infty) \times N$. 
$\beta_{\lambda, \mu}$ is invariant under the 
Hopf action 
of $H(t)$ for $t\in \frac{2\pi}{3}\Z$. 
On the end $(T, \infty) \times N$ naturally 
$\beta_{\lambda, \mu}$ admits more symmetries, namely,  
it is invariant under the translations 
in $\tau$-direction and also under the Hopf flow 
$H(t)$ for any $t\in \R$.   
}
\end{cor}
Apart from the main result Theorem~\ref{MainTheorem} 
of this article, this result might be of an independent 
interest.  In the final section, 
we will make a brief discussion 
on the generalization of this result, 
namely, the (non-)existence of 
an end-periodic symplectic structure 
on Stein or globally convex symplectic manifolds.  
In the rest of this section, 
we prove the above theorem.  
\begin{lem}\label{seed}
{\rm\quad
There exists a closed 2-form $\kappa$ on 
${\C}P^2 \setminus \E$ which restricts to 
$ dy\wedge \zeta_{N'}$ on the product end. 
}
\end{lem}

{\it Proof} of Lemma~\ref{seed}. \quad 
As $ dy\wedge \zeta_{N'}$ is closed, it defines 
a de Rham cohomology class 
$[dy\wedge \zeta_{N'}]\in H^2(N'
)\cong \R^2$. 
Let us look at the Meyer-Vietoris exact sequence 
for the cohomologies of 
${\C}P^2 = \ol{U} \cup ({\C}P^2\setminus U$). 
It is easy to see that the inclusion 
$N' \hookrightarrow  \ol{U}$ 
induces a trivial map $0: H^2(U) 
{\to} H^2(N') $. 
Therefore the fact $H^3({\C}P^2)=0$ and 
the long exact sequence tell us that the inclusion  
to the other side induces a surjective homomorphism 
$H^2({\C}P^2\setminus U) \twoheadrightarrow H^2(N') $.  
This also implies that 
the closed 2-form 
$dy\wedge \zeta_{N'}$ 
on the product end 
extends to the whole ${\C}P^2 \setminus \E$ as 
a closed 2-form $\kappa$.  \hfill $\square$\ref{seed}.
\begin{rem}\label{SimpleEllipticRem}
{\rm\quad 
The origin of ${\bf F}_0$ is an isolated singularity 
of {\it simple elliptic} type.  
Up to 3-fold branched covering $U$ is 
orientation-reversing diffeomorphic 
to the minimal blowing up resolution of ${\bf F}_0$.  
On the other hand, ${\C}P^2\setminus\E$ is up to 
3-fold covering biholomorphic 
to the Milnor fibre. 
The above lemma reflects 
the fact that the resolution and the Milnor fibre 
are quite different to each other.  
Such a phenomenon does not happen for 
{\it simple} singularities.  
For singularity theory, see for example 
\cite{D}.  
}
\end{rem}

Let us proceed to construct 
an end-periodic symplectic form 
on  ${\C}P^2 \setminus \E$. 
First take a positive constant $\mu$ small enough so that 
$\beta' + {\mu}\,{dy}\wedge\zeta_{N'}$ 
is still a symplectic form.  
On the product end $\{\tau\geq 2\pi\}$, 
from Corollary 4.2 we know   
$\beta'= d(e^\tau\zeta_{N'})
=e^\tau d\tau\wedge\zeta_{N'} 
+e^\tau\frac{3}{2\pi}dx\wedge dy$.  
This implies 
 $\beta'\wedge dy\wedge\zeta_{N'}=0$ 
and hence we have 
 $(\beta' + \mu \, dy\wedge\zeta_{N'})^2 = {\beta'}^2$ 
on the product end. 
Therefore if we choose $\mu$ small enough, we can assure 
that even on the compact core ${\C}P^2 \setminus U$, 
the closed 2-form $\beta' + \mu\, dy\wedge\zeta_{N'}$ 
is non-degenerate. We fix such $\mu$. 

Next take constants 
$2\pi<T_0<T_1<T_2<T_3=T$ 
and 
non-negative smooth functions 
$k(\tau)$ and $l(\tau)$ of $\tau$ 
on $[T_0, \infty)$ satisfying the following conditions:   
$$
\begin{array}{rrl}
k(\tau)=e^\tau, 
&\quad 
l(\tau)\equiv0\: :
&
\quad  T_0\leq\tau\leq T_1 , 
\\
k'(\tau) > 0\,,
&
l(\tau) >0\: :
&
\quad  T_1\leq\tau < T_2 ,
\\
k(\tau) >0\,,
&
l(\tau) \equiv \lambda\: :
&
\quad  T_2\leq\tau\leq T_3 ,
\\
k(\tau) \equiv0\,,
&
l(\tau) \equiv \lambda\: :
&
\quad  T_3\leq\tau .
\end{array}
$$
This is done as follows.  
First choose such a smooth function $k$. 
Then take a constant $\lambda>0$ 
satisfying 
$\max\{-\frac{3k'(\tau)k(\tau)}{4\mu\pi}\, 
;\, T_2 \leq \tau \leq T_3 \} < \lambda$.  
Then it is easy to find a smooth function $l(\tau)$ 
which satisfies all of the above conditions. 

Now we are ready to construct 
an end-periodic symplectic form. 
First modify $\beta'$ on the product end.  
We can define $\beta'_\sharp$ as
\vspace{-5pt}
$$
\beta'_\sharp = \left\{ 
\begin{array}{ll}
\beta' \quad 
& 
\mathrm{on}\quad {\C}P^2\setminus U \, ,
\\ 
d(k(\tau)\zeta_{N'})
+
l(\tau)d\tau\wedge dx
&
\mathrm{on}\quad [T_0, \infty) \times N'
\end{array}
\right. 
\vspace{-5pt}
$$
because two expressions of $\beta'_\sharp$ coincide 
with each other on $[T_0, T_1]\times N'$.  
Finally we put 
$\beta'_{\lambda, \mu}=\beta'_\sharp + \mu \kappa$.  
This is the desired symplectic form 
on ${\C}P^2\setminus \E$ 
because of the following reasons. 
First of all, $\beta'_{\lambda, \mu}$ is closed and 
coincides with 
$\lambda\,d\tau\wedge dx + \mu\,dy\wedge \zeta_{N'}$ 
on $[T_3, \infty)\times N'$ 
and with $\beta' + \mu\,dy\wedge\zeta_{N'}$ 
on ${\C}P^2\setminus U$. 
Therefore it is non-degenerate 
on ${\C}P^2\setminus U$ as already remarked above. 
On the product end, 
as 
$d(k(\tau)\zeta_{N'})$ 
and 
$l(\tau)d\tau\wedge dx + \mu\, dy\wedge\zeta_{N'}$ 
do not interact at all under the exterior product,  
we have  
$$
(\beta'_{\lambda, \mu})^2 
=
\left( 
\frac{3k'(\tau)k(\tau)}{2\pi}
+ 2l(\tau)\mu
\right)
d\tau\wedge dx \wedge dy \wedge \zeta_{N'} .
$$
Therefore 
$\beta'_{\lambda, \mu}$ 
is non-degenerate on the product end as well. 
\hfill $\square$\ref{end-periodicThm}.

\section{\large $\widetilde{E}_7$ 
and $\widetilde{E}_8$}\label{SimpleElliptic}

Among simple elliptic singularities, 
the following three deformation classes 
$\widetilde{E}_l$ ($l=6,7,8$) 
are known to be realized as hypersurface 
singularities and their links are ismorphic 
to $\nil{-3}$,  $\nil{-2}$, and to $\nil{-1}$ 
respectively. They are defined by the 
following polynomials. 
\begin{eqnarray*}
f_{\widetilde{E}_6} &=& Z_0^3 + Z_1^3 + Z_2^3 
\,\,\,(+ \lambda Z_0Z_1Z_2) 
\\
f_{\widetilde{E}_7} &=& Z_0^4 + Z_1^4 + Z_2^2 
\,\,\,(+ \lambda Z_0Z_1Z_2) 
\\
f_{\widetilde{E}_8} &=& Z_0^6 + Z_1^3 + Z_2^2 
\,\,\,(+ \lambda Z_0Z_1Z_2) 
\end{eqnarray*}
As the smooth topology of these objects does not 
depend on the choice of the constant $\lambda$ 
(except for finitely many exceptional values), 
in this paper we ignore it and take it to be $0$.  
Each of our constructions in this paper 
for the Fermat cubic, \IE 
the $\widetilde{E}_6$ polynomial, also works in  
the other two cases.  
In this section we verify this fact,  
by briefly reviewing the topology 
of these singularities.  
For basic facts about hypersurface singularities,  
the readers may refer to 
Milnor's seminal text book \cite{M}.  
Also Dimca's book \cite{D} provides more detailed 
and modern basic informations. 
\vspace{10pt}

The notations in \S\ref{Lawson} are 
used and interpreted in parallel 
or slightly modified 
meanings according to the context,  
unless otherwise specified.  
For 
$f=f_{\widetilde{E}_l}$ 
($l=7,8$) 
the origin is an isolated and 
in fact unique singularity of the hypersurfaces 
${\bf F}_0$.  
Instead of scalar multiplication, 
we define the weighted homogeneous action 
of $\lambda\in\C^\times$ on $\C^3$ 
by 
$\lambda\cdot(Z_0, Z_1, Z_2)
=(\lambda^{w_0}Z_0, \lambda^{w_1}Z_1, 
\lambda^{w_2}Z_2)$ 
where the weight vector 
${\bf{w}}=(w_0, w_1, w_2)$ takes value  
$(2,1,1)$ \RESP{$(3,2,1)$} 
for 
$\widetilde{E}_l$  ($l=7,8$). 
By this action  we have 
$\lambda\cdot {\bf F}_w={\bf F}_{\lambda^4w}$ 
\RESP{$\lambda\cdot {\bf F}_w={\bf F}_{\lambda^6w}$} 
for 
$l=7$ 
\RESP{$l=8$}.  
The weighted homogeneous action 
by positive real numbers 
$\lambda\in\R_+$ plays the role 
of the euclidean homotheties 
in the Feramt cubic case. 
The action by unit complex numbers 
$\lambda=e^{it}$ ($t\in\R$) restricts to 
the action on $S^5$ and is again denoted 
by $H(t)$ and called the 
{\it weighted Hopf} action or flow.  
The quotient space
$\WP$ 
of this weighted Hopf 
action 
is called 
the {\it weighted projective space}, 
which is a complex analytic orbifold.  
The quotient map $h:S^5 \to \WP$ 
is called the {\it weighted Hopf fibration}, 
which is a Siefert fibration.  
We take $\C P^2=\{[X_0:X_1:X_2]\}$ 
as a quotient of $\WP$ as follows.  
Define a map 
$\Phi:\cC^3 \to \C P^2$ 
as 
$\Phi:(Z_0,Z_1, Z_2) \mapsto 
[X_0:X_1:X_2] 
=
[Z_0^{d_0}:Z_1^{d_1}:Z_2^{d_2}]$ 
where 
$ 
(d_0, d_1, d_2)
=
(1,2,2)
$
\RESP{$(2,3,6)$} 
for $\widetilde{E}_7$ 
\RESP{$\widetilde{E}_8$}. 
Then 
$\Phi$ 
factors into $\Phi\vert_{S^5}=\Psi\circ h$ 
for some $\Psi:\WP \to \C P^2$.  
The homogeneous equations 
\vspace{-5pt}
\begin{eqnarray*}
g_{\widetilde{E}_7} &=& X_0^2 + X_1^2 + X_2^2 \,\, 
=\,\, 0
\\ 
g_{\widetilde{E}_8} &=& X_0 + X_1\, + X_2 \,\,\,
=\,\,0\, 
\end{eqnarray*}
on $\C P^2$ rewrites 
$f_{\widetilde{E}_l}=0$ 
as 
$g_{\widetilde{E}_l}\circ\Phi=0$ ($l=7,8$).

The first important fact to notice is 
that the open set of $S^5$ consisting of all regular 
orbits of the weighted Hopf flow $H(t)$ 
contains $N={\bf F}_0\cap S^5$. Therefore the orbit space 
$E_{(l)}=
N/_H$
is a non-singular 
holomorphic curve which sits in the regular part of 
$\WP$.

`$g_{\widetilde{E}_7}=0$' defines a non-singular 
projective curve of degree 2 and 
`$g_{\widetilde{E}_8}=0$' a projetive line. 
Both of them are biholomorphic to $\C{P}^1$.  
Comparing $\Phi$ and $h\vert_{E_{(l)}}$, 
we easily see 
that $\Psi\vert_{E_{(7)}}:E_{(7)}\to 
\{X_0^2 + X_1^2 + X_2^2=0\}$ 
is a 2-fold branched covering 
over the rational curve with 
4 branched points 
$
\{X_1=0\,\, \mathrm{or }\, X_2=0\}
\cap
\{X_0^2 + X_1^2 + X_2^2=0\}
$ 
like the Weierstrass $\wp$ function and 
$E_{(7)}$ is seen to be an elliptic curve. 
In the case of $\widetilde{E}_{(8)}$,  
$\Psi\vert_{E_{(8)}}:E_{(8)}
\to 
\{X_0 + X_1\, + X_2=0\}$ 
is a 6-fold branched covering, 
branching over 3 points 
$\{X_0=0\}$, $\{X_1=0\}$, and $\{X_2=0\}$ 
with branch indices 2, 3, and 6 
respectively. From this we also see that 
$E_{(8)}$ is an elliptic curve. 

Similarly it is easy to see that 
the self-intersection 
(the $c_1$ of the normal bundle) 
of $E_{(7)}$ 
\RESP{$E_{(8)}$} 
in $\WP$ 
is 8 
\RESP{6}  
and that 
the 
$c_1$ of the weighted Hopf fibrations over 
$E_{(7)}$ 
\RESP{$E_{(8)}$} 
is 
$-2$ \RESP{$-1$}. 

Like in the case of $\widetilde{E}_{(6)}$, 
in both of the other two cases the weighted projection 
$\cC^3 \to S^5=\cC^3/\R_+$ 
by positive real numbers 
restricts to a diffeomorphism from
${\bf F}_{re^{i\theta}}$ 
to the Milnor fibre 
$L_\theta$.  
$h\vert_{L_\theta}:L_\theta \to 
\WP\setminus \widetilde{E}_{(l)}$ 
is a branched covering, 
but the number of branched points is finite 
and around the ends it is a 4-fold \RESP{6-fold} 
regular covering for $l=7$ \RESP{$l=8$}. 

We also remark here that the link $N$ has a 
product type tubular neighbourhood $W$ in $\WP$ 
because $f$ gives the trivialization. 
The boundary $\partial W$ is a Kodaira-Thurston 
nil-manifold and 
$\cF_0$ can be considered as its cyclic covering.  
\vspace{10pt}

Now let us verify that 
our constructions are transplanted to 
the cases of $\widetilde{E}_7$ and $\widetilde{E}_8$. 
From the descriptions of the link $N$, 
the Milnor fibres $L_\theta$, and of ${\bf F}_1$, 
the contents in Section 1 and 2 are recovered. 
The fact that 
the weighted Hopf flow preserves 
the standard sphere $S^5(\rho)$ of radius $\rho$ 
and the standard symplectic form $\beta^*$ 
implies that 
the product coordinates introduced in Section 3 
can play the same role as in the case of $\widetilde{E}_6$ 
and the arguments in Section 4 are valid without modifications. 

As to the results in Section 5, 
once a parallel result to Lemma \ref{seed} 
is verified, then the manipulations of 
differential forms on the product end 
holds without major modifications. 
Together with the commutative 
diagram below, 
the fact 
that the rational (or real) cohomology 
of $\WP$ is isomorphic to that of $\C{P}^2$ 
(see \EG \cite{D}) 
tells 
that a parallel to Lemma \ref{seed} 
holds.  
$$
\begin{array}{ccccc}
(T, \infty)\times N \,\,\,& \cong 
& 
\mathrm{end \,\, of\,\, }\cF_0 
\cong 
\mathrm{end \,\, of\,\, }{\bf F}_1 
& \hookrightarrow & {\bf F}_1\cong L_0
\vspace*{2pt}
\\
\,\,\,\,\,\,\,\,
{\big\downarrow} & & \big\downarrow & & 
\,\,\,\,\,\,\,
\big\downarrow{\scriptstyle{h\vert_{L_0}}}
\vspace*{3pt}
\\
(T, \infty)\times \partial U & \cong & 
\mathrm{end \,\, of\,\, }{\WP\setminus E_{(l)}} 
& 
\hookrightarrow 
&
\WP\setminus E_{(l)}
\end{array}
$$
The left and the middle vertical arrows 
are regular coverings and the right one 
is a branched covering.

\section{\large Concluding Remarks}\label{concluding}

To close the present article, we make some comments 
and raise some questions related to our construction.  

\subsection{End-periodic symplectic structures on 
Stein or globally convex symplectic manifolds}

The construction of leafwise symplectic structure 
in this paper seems to stand on 
an extremely rare intersection of fortunes. 

Besides the foliations of codimension one, 
as is mentioned in the previous section, 
the existence of end-periodic 
symplectic structures on Stein or globally convex 
symplectic manifolds might be of an independent 
interest.  
However, the possibility of such cases 
seems to be still limited.  
In the final section we discuss on such problems. 
\vspace{5pt}

\begin{ex}\label{C}
{\rm \quad 
The Stein manifold $\C$ 
(or the upper half plane $\Bbb H$) 
carries an end-periodic symplectic form. 
}
\end{ex}
This example is in many senses trivial, 
because, 
first of all the fact itself is trivial.  
Especially we do not have to change the 
symplectic form.  
Also, as this Stein 
manifold is not really convex, we should say  
this is a meaningless example.  
The convexity 
of symplectic structures must be discussed 
on manifolds of dimension $\geq 4$ (see \cite{EG}).  
However, this example still exhibits a 
clear contrast to the following example. 

\begin{ex}\label{C^n}
{\rm \quad 
The Stein manifold $\C^n$ ($n\geq2$) 
does not admit an end-periodic symplectic 
structure, basically because 
$S^1\times S^{2n-1}$ does not admit 
such structures.   
}
\end{ex}
As the non-existence 
this example is 
generalized to many cases. 
\vspace{5pt}

For the case of symplectic dimension 4, 
the recent result by Friedl and Vidussi, 
which has been known as Taubes' conjecture, 
provides a strong constraint. 

\begin{thm}{\rm
(Taubes' conjecture, Friedl-Vidussi \cite{FV})
\quad 
For a closed 3-manifold $M$, 
the 4-manifold $W^4=S^1\times M^3$ admits 
a symplectic structure 
if and only if $M$ fibers over the circle. 
}
\end{thm}

Now in order to make 
the implication of this theorem clearer, 
let us take the following definition.  
\begin{df}
{\rm \quad 
Assume that an open $2n$-manifold $W$ 
has an end which is diffeomorphic to 
$\R_+\times M^{2n-1}$ for some closed oriented 
manifold $M$.   
An {\it end-periodic symplectic structure} 
on $W$ is a symplectic structure on $W$ 
whose restriction to the end is invariant under 
the action of non-negative integers $\N_0$ 
where $m\in\N_0$ acts on $\R_+\times M$ 
as $(t, x)\mapsto (t+m, \varphi^m(x))$ 
for some fixed monodromy diffeomorphism 
$\varphi:M\to M$.  
}
\end{df}
It follows directly from the definition that 
the mapping cylinder $M_\varphi$ admits 
a symplectic structure. 
If the monodromy belongs to a mapping class of 
finite order, $S^1\times M$ also admits a 
symplectic structure. 

Above theorem due to Friedl and Vidussi 
tells that in the trivial or finite monodromy 
case, the $M^3$ must fibres over the circle. 
Essentially the same construction 
in the case of the Kodaira-Thurston 
nil-manifold gives  symplectic structures  
on such closed 4-manifolds. 
The virtue of their theorem is of course in 
the converse implication. 

\begin{ex}{\rm (Surfaces of higher degree)\quad 
Instead of taking the Fermat cubic 
surface as in the present article, 
if we take, \EG a Fermat type quartic surface, 
then the end is diffeomorphic to 
$\R_+\times M^3$, 
where $M^3$ is an $S^1$-bundle over 
the closed oriented surface  $\Sigma_3$  
of genus 3 with euler class 4. 
As this 3-manifold apparently does not fiber 
over the circle, 
the Fermat quartic surface does not 
admit an end-periodic symplectic structure with 
at most finite monodromy. 
The same applies 
to the Brieskorn type hyperbolic singularities 
$\{Z_0^p+Z_1^q+Z_2^r=0\}$ 
with 
$1/p + 1/q + 1/r < 1$.  
}
\vspace{5pt}
\end{ex}
\begin{ex}{\rm (Cusp singularities)\quad 
For $1/p + 1/q + 1/r <1$, 
the polynomials equation 
$Z_0^p+Z_1^q+Z_2^r+Z_0Z_1Z_2  =0$ 
defines so called the {\it cusp} singularity at the origin. 
(A cubic term is added to the above Brieskorn type polynomial.)
The link is a solv-3-manifold, which are explained in the next example. 
This seems to be only the possible case other than what 
we treated in this paper, where the Milnor fibration 
is modified into a codimension 1 leafwise symplectic foliation. 
In fact it turns out to be possible.  
However, as it requires totally different arguments, it is 
discussed in a forthcoming paper.  
\vspace{5pt}}
\end{ex}

If we extend our scope from Stein manifolds to 
globally convex symplectic manifolds (see \cite{EG} 
for details of this notion) 
we find one more example for the existence 
of end-periodic symplectic sturcture, 
which is only slightly less trivial than Example~\ref{C}.  
\begin{ex}\label{solv}
{\rm (Solvable manifold)\quad 
Let $M^3=$~{\it Solv} be a 3-dimensional 
compact solv-manifold, 
namely, 
the mapping cylinder 
of a hyperbolic automorphism of $T^2$. 

Then it carries an algebraic Anosov flow, 
whose strong (un)stable direction 
corresponds  to an eigenvalue.  
Then from \cite{Mi} we know that $\R\times M$ 
admits a globally convex symplectic structure. 
As it has a disconnected end, 
it is not Stein. 
On the other hand, apparently  $\R\times M$ admits 
an end-periodic symplectic structure because 
$S^1\times M$ is symplectic.  
}
\end{ex}
Some of higher genus surface bundles over the circle 
with pseudo-Anosov monodromy 
admit Anosov flows.  See \cite{Go} and \cite{FH} for 
Anosov flow on hyperbolic 3-manifolds. 
Of course in this formulation, 
Thurston's conjecture (now it is a theorem) is involved, 
which asserts that 
any closed hyperbolic 3-manifold admits 
a finite covering which fibers over the circle. 
The above example is extended to those. 

These examples are again disappointing 
because 
the compact core part has no more topology than 
the end.  
\vspace{5pt}

Apart from trivial constructions 
like taking products of the Fermat cubic surfaces,  
Examples \ref{C} and \ref{solv} are the only 
known other such manifolds that is  
convex and at the same time 
admits an end-periodic symplectic structure. 
Also remark that 
the construction of differential forms 
in Section~\ref{end-periodic} might have 
some similarities only on 
manifolds of dimension $4k$ ($k\in\N$). 
\vspace{10pt}

\subsection{Foliations on spheres}

Meersseman and Verjovsky proved in \cite{MV2} 
that Lawson's foliation 
did not admit a leafwise complex structure. 
The tube component is obstructed to admit 
a leafwise complex structure. 
Contrary to the case of 
leafwise symplectic structures, 
it is easy to show that the Fermat cubic component 
admits a leafwise complex structure but 
for the tube comonent it is the other way round. 
At present, it still seems to the author that 
the existence of a foliation 
of codimension-one on $S^5$ 
with smooth leafwise complex structure 
is still a totally open problem.  
\vspace{5pt}

Also the existence of codimension-one foliation 
with leafwise symplectic structures 
on higher odd dimensional spheres is unknown as well.  
\vspace{5pt}

A recent result by G. Meigniez (\cite{Me}) 
claims that in dimension $\geq 4$, 
once a smooth foliation of codimension one 
exists on a closed manifold, it can be modified 
into a minimal one, namely, 
one with every leaf dense. 
Especially it follows that 
in higher dimensions there is no more 
direct analogue of 
Novikov's theorem, that is, 
for any foliation of codimension 1 
on $S^3$ there exists a compact leaf.  
It is also mentioned by Meigniez that 
under the presence of some geometric structures, 
like leafwise symplectic structures,  
it might be of some interest to ask whether 
some similar statement to Novikov's theorem 
holds or not.

\begin{flushright}
Yoshihiko MITSUMATSU 
\vspace{3pt}
\\
{\small Department of Mathematics, Chuo University
\\
1-13-27 Kasuga Bunkyo-ku, Tokyo, 112-8551, Japan
\\
e-mail : yoshi@math.chuo-u.ac.jp
\vspace{3pt}
\\
}
\end{flushright}

\end{document}